\documentclass[12pt]{article}

\usepackage{amsmath, amssymb,amsthm,amsbsy}
\usepackage{graphicx}
\usepackage{psfrag}

\numberwithin{equation}{section}

\newtheorem{theorem}{Theorem}[section]
\newtheorem{lemma}[theorem]{Lemma}

\newtheorem{corollary}[theorem]{Corollary}

\theoremstyle{definition}
\newtheorem{example}[theorem]{Example}

\newtheorem{remark}[theorem]{Remark}

\newtheorem*{acks}{Acknowledgements}

\theoremstyle{remark}

\newcommand{\refT}[1]{Theorem~\ref{#1}}
\newcommand{\refC}[1]{Corollary~\ref{#1}}
\newcommand{\refL}[1]{Lemma~\ref{#1}}

\newcommand{\refS}[1]{Section~\ref{#1}}

\newcommand{\refE}[1]{Example~\ref{#1}}

\newenvironment{romenumerate}[1][0pt]{
\addtolength{\leftmargini}{#1}\begin{enumerate}
 }{\end{enumerate}}

\newcounter{oldenumi}
{\setcounter{oldenumi}{\value{enumi}}
\begin{romenumerate} \setcounter{enumi}{\value{oldenumi}}}
{\end{romenumerate}}

\newcounter{thmenumerate}

\newcounter{xenumerate}   

\newcommand\nfrac{\tfrac} 

\def\urltilde{\kern -.15em\lower .7ex\hbox{\~{}}\kern .04em}

\begingroup
  \divide\count255 by 60
  \multiply\count255 by -60
  \advance\count255 by \time
  \ifnum \count255 < 10 \xdef\klockan{\the\count1.0\the\count255}
  \else\xdef\klockan{\the\count1.\the\count255}\fi
\endgroup

\newcommand\set[1]{\ensuremath{\{#1\}}}
\newcommand\bigset[1]{\ensuremath{\bigl\{#1\bigr\}}}
\newcommand\Bigset[1]{\ensuremath{\Bigl\{#1\Bigr\}}}

\newcommand\bigpar[1]{\bigl(#1\bigr)}
\newcommand\Bigpar[1]{\Bigl(#1\Bigr)}

\newcommand\lrpar[1]{\left(#1\right)}

\newcommand\Bigabs[1]{\Bigl|#1\Bigr|}

\def\rompar(#1){\textup(#1\textup)}    

\newcommand\bigparfrac[2]{\bigpar{\frac{#1}{#2}}}
\newcommand\Bigparfrac[2]{\Bigpar{\frac{#1}{#2}}}

\def\xexp(#1){e^{#1}}

\newcommand\ntoo{\ensuremath{{n\to\infty}}}

\newcommand\ie{i.e.\spacefactor=1000}

\newcommand{\tend}{\longrightarrow}
\newcommand\dto{\overset{\mathrm{d}}{\tend}}

\newcommand\bbR{\mathbb R}

\newcommand\bbZ{\mathbb Z}
\newcommand\bbQ{\mathbb Q}

\newcounter{CC}
\newcommand{\CC}{\stepcounter{CC}\CCx} 
\newcommand{\CCx}{C_{\arabic{CC}}}     
\newcounter{cc}
\newcommand{\cc}{\stepcounter{cc}\ccx} 
\newcommand{\ccx}{c_{\arabic{cc}}}     
\newcommand{\ccdef}[1]{\xdef#1{\ccx}}     

\newcommand\E{\operatorname{\mathbb E{}}}
\renewcommand\P{\operatorname{\mathbb P{}}}

\newcommand\rank{\operatorname{rank}}
\newcommand\Vol{\operatorname{det}}

\newcommand\ga{\alpha}
\newcommand\gb{\beta}
\newcommand\gd{\delta}

\newcommand\gG{\Gamma}

\newcommand\gl{\lambda}

\newcommand\gs{\sigma}

\newcommand\eps{\varepsilon}

\newcommand\cL{{\mathcal L}}

\newcommand\ett[1]{\boldsymbol1[#1]}
\newcommand\etta{\boldsymbol1}
\def\[#1]{[\![#1]\!]}

\newcommand\qq{^{1/2}}
\newcommand\qqq{^{1/3}}

\newcommand\qqqw{^{-1/3}}
\newcommand\qqw{^{-1/2}}
\newcommand\qw{^{-1}}
\newcommand\qww{^{-2}}

\renewcommand{\=}{:=}

\newcommand\oi{[0,1]}

\newcommand\dd{\,\textup{d}}

\newcommand\half{\frac12}
\newcommand\thalf{\tfrac12}

\newcommand\innprod[1]{\langle #1\rangle}

\newcommand{\inre}{^\circ}
\newcommand{\tqn}{\widetilde Q_n}
\newcommand{\sqrtn}{\sqrt{n}}
\newcommand{\lln}{\cL+\ell_n}
\newcommand{\llnx}{(\lln)}
\newcommand{\K}{K}
\newcommand{\bbrn}{\bbR^N}
\newcommand{\bbzn}{\bbZ^N}

\newcommand{\xxm}{x_1,\dots,x_m}
\newcommand{\xxg}{x_1,\dots,x_g}

\newcommand{\ttm}{t_1,\dots,t_m}
\newcommand{\ort}{^\perp}

\newcommand{\lnx}[1]{{L}_n(#1)}
\newcommand{\lng}{\lnx{G}}

\newcommand{\gv}{{g}}

\newcommand{\pmx}{perfect matching}
\newcommand{\fpm}{fractional perfect matching}
\newcommand{\elle}{\ell_e}
\newcommand{\hA}{\widehat{A}}
\newcommand{\T}{^T}
\newcommand{\clg}{\cL^{(1)}_G}
\newcommand{\clv}{\cL^{(2)}_G}
\newcommand{\DG}{D_G}
\newcommand{\rv}{|_V}

\newcommand{\Gram}[2]{(\innprod{#1_i,#1_j})_{i,j=1}^{#2}}
\newcommand{\detGram}{\det\Gram}
\newcommand{\ief}{_{ief}}
\newcommand{\aax}{\vec A}
\newcommand{\aaxg}{\aax_G}
\newcommand{\ax}{\vec a}
\newcommand{\dsx}{\tilde D_s}

\newcommand{\maple}{\texttt{Maple}}

\newcommand\REM[1]{{\raggedright\texttt{[#1]}\par\marginal{XXX}}}

\begin{document}

\title{On the number of perfect matchings in random lifts }

\author{
Catherine Greenhill \\
\small School of Mathematics and Statistics\\[-0.8ex]
\small University of New South Wales\\[-0.8ex]
\small Sydney, Australia 2052\\[-0.3ex]
\small\texttt{csg@unsw.edu.au}
\and
Svante Janson \\
\small Department of Mathematics\\[-0.8ex]
\small Uppsala University\\[-0.8ex]
\small PO Box 480 \\[-0.8ex]
\small S-751 06 Uppsala, Sweden\\[-0.3ex]
\small \texttt{svante.janson@math.uu.se}
\and
Andrzej Ruci{\' n}ski\thanks{Supported by grant N201 036 32/2546}\\
\small Department of Discrete Mathematics\\[-0.8ex]
\small Adam Mickiewicz University\\[-0.8ex]
\small Pozna{\' n}, Poland 61-614\\[-0.3ex]
\small \texttt{rucinski@amu.edu.pl}
}

\date{2 July 2009; revised 5 November 2009}

\maketitle

\begin{abstract}
Let $G$ be a fixed connected multigraph with no loops.
A random $n$-lift of $G$
is obtained by replacing each vertex of $G$ by a set 
of $n$ vertices 
(where these sets are pairwise disjoint)
and replacing each edge by a randomly chosen perfect matching
between the $n$-sets corresponding to the endpoints of
the edge.
Let $X_G$ be the number of perfect matchings in a random lift
of $G$.
We study the distribution of $X_G$ in the limit as $n$ tends to
infinity, using the small subgraph conditioning method.

We present several results including an asymptotic
formula for the expectation of $X_G$ when $G$ is $d$-regular,
$d\geq 3$.  
The interaction of perfect matchings with short cycles
in random lifts of regular multigraphs is
also analysed.  Partial calculations are performed for the
second
moment of $X_G$, with full details given for two example 
multigraphs, including the complete
graph $K_4$.

To assist in our calculations we provide a theorem for estimating
a summation over multiple dimensions using Laplace's method.
This result is phrased as a summation over lattice points,
and may prove useful in future applications.
\end{abstract}


\section{Introduction}\label{intro}

Throughout, let $G$ be a fixed connected multigraph with 
$\gv$ vertices and no loops.
For simplicity we assume that $V(G)=[\gv]\=\set{1,\dots,\gv}$.
A random $n$-lift of $G$ is a random graph on the vertex set
$V_1\cup V_2\cup \cdots \cup V_\gv$, where each $V_i$ is a set
of $n$ vertices and these sets are pairwise disjoint, obtained by
placing a uniformly chosen random perfect
matching between $V_i$ and $V_j$, independently for each
edge $e=ij$ of $G$.  
Denote the resulting random graph by $\lng$.
The perfect matching corresponding to the
edge $e$ of $G$ is called the \emph{fiber} corresponding to $e$,
which we denote by $F_e$. 
Note that the degree of $v\in V_i$ in $\lng$ is equal to
the degree $d_G(i)$ of vertex $i$ in $G$.
In particular, if $G$ is $d$-regular, then so is $\lng$.
We are interested in asymptotics as $n$ tends to infinity.

This model of sparse random graphs was introduced and studied in 
a series of papers by Amit, Linial, Matou{\v s}ek, and 
Rozenman~\cite{AL1,AL2,ALM,LR}. Linial and Rozenman~\cite{LR} studied
the existence of a perfect matching in $\lng$ and described a large class 
of graphs $G$ for which $\lng$ a.a.s.\  contains a perfect matching 
(for $n$ even, at least). This class contains all regular graphs and,
in turn, is contained in the class of graphs having a fractional perfect
matching (see Section~\ref{exppm} for a definition). 
Observe that if $G$ has a perfect matching then every lift of $G$ has 
at least one perfect matching.

In this paper we study the number of perfect matchings
in $\lng$ in the limit as $n$ tends to infinity, 
where $G$ is a graph with a fractional perfect matching.
To do this we use the \emph{small subgraph conditioning method},
which provides a concentration result based on the 
second moment method conditioned on the number of small cycles. 
For a concise description of the 
method, see \cite[Theorems 9.12 and 9.13]{JLR}.

Let $X_G$ be the number of perfect matchings in $\lng$.
To apply the small subgraph conditioning method, asymptotic
expressions for $\E X_G$ and $\E(X_G^2)$ must be found.   Then
the limit of the ratio $\E(X_G^2)/(\E X_G)^2$ is compared against
a quantity which
depends upon the interaction of perfect matchings and short cycles
in $\lng$.  

In Sections~\ref{exppm} and~\ref{sm} we write the first and second
moments of $X_G$ as multiple sums of some explicit terms, and then estimate
the sums by Laplace's method. This is a standard method for similar moment 
estimates, 
and in particular, it has been used in several papers on random regular graphs. 
(See for example \cite[Chapter 9]{JLR} and the references given there.)
However, in the present paper, each summation is over an index set of rather 
high dimension with a number of side conditions on the indices, while in many 
previous applications the summations are only over one or two variables.
To assist with these calculations, we present 
a general theorem (\refT{TA2}) that
encapsulates Laplace's method for a general situation, with sums over a lattice 
in a subspace of $\bbR^N$. We do this both because we think that it clarifies 
the argument in the present work, and
because we hope that it might be useful in future applications.
The necessary terminology and notation is introduced in Section~\ref{notation},
where Theorem~\ref{TA2} is stated.  The proof of Theorem~\ref{TA2}
can be found in Section~\ref{appa}.

Using this machinery we prove an asymptotic formula
for $\E X_G$ for any connected regular multigraph $G$ with degree at 
least three (see Theorem~\ref{expectation}). 
However, two difficulties (one algebraic and one analytic)
have prevented us from obtaining
an asymptotic formula for $\E(X_G^2)$ in the same generality,
though we have partial results in
Theorem~\ref{conditional} and Lemma~\ref{Llattice}.
We illustrate these results by calculating $\E(X_G^2)$ for
two multigraphs: specifically, for the complete graph 
$K_4$ and for the multigraph consisting of two vertices
and three parallel edges, which we denote by $K_2^3$.
These calculations were performed with the aid of \maple.
(A file containing the \maple\ commands is available 
from~\cite{maplecalcs}.)

In Section~\ref{concen} we prove the necessary results relating
to short cycles in random lifts (Lemmas~\ref{lambdak}, \ref{muk}
and Corollary~\ref{Cw}). 
As corollaries, using~\cite[Theorem 9.12]{JLR} we obtain
a concentration result for $X_G$ in our two illustrative examples
(see Corollaries~\ref{mainK4} and \ref{mainbanana}).

One of the most interesting questions on random lifts is the problem of 
existence of a Hamilton cycle.   
There is a conjecture (attributed to Linial)
that a random lift of $K_4$ is a.a.s.\ hamiltonian.  Indeed,
we believe that a.a.s.\ $\lnx{G}$ is hamiltonian for all connected
$d$-regular loop-free multigraphs $G$ with $d\geq 3$.
(This is known to be true when $G$ is a multigraph with 
exactly two
vertices and at least three edges: see Remark~\ref{para} below.) 
Burgin, Chebolu, Cooper and Frieze~\cite{BCCF} showed that 
a.a.s.\ $\lnx{K_g}$ is hamiltonian when $g$ is large enough
(see also \cite{CF} for the directed case).  The arguments in
\cite{BCCF} are combinatorial and utilize the
celebrated  idea of P\'osa.  For small $g$, we feel that
the small subgraph conditioning method may be a fruitful line
of attack, as it has been very successful for
studying Hamilton cycles in random regular graphs 
(Robinson and Wormald \cite{RW92,RW94}, see also \cite[Chapter 9]{JLR}). 
This remains an open problem.

\begin{remark}\label{para}
We allow the multigraph $G$ to have multiple edges. The simplest case is when 
$G$ consists of only two vertices, with $d$ parallel edges between them. 
The random lift $\lng$ then is a random bipartite (multi)graph obtained 
by taking 
the union of $d$ independent random matchings between two sets of
$n$ vertices each. Such sums have been studied in \cite{MRRW}, where they were 
shown to be contiguous to random bipartite $d$-regular (multi)graphs. 
The latter, in turn, is known to be a.a.s.\ hamiltonian (see \cite{RW-bip} for 
a standard, second moment method proof).  Hence for this small multigraph
$G$ with $d\geq 3$, the random lift $\lng$ is a.a.s.\ hamiltonian too.
\end{remark}

\begin{remark} 
Random lifts of multigraphs with loops can also be formed.
As in~\cite{AL1}, the fiber corresponding to a loop 
is given by the $n$
  edges $i\gs(i)$ for a random permutation $\gs$ of $[n]$.  
This is a
  random 2-regular (multi)graph, denoted by $\mathbb{P}(n)$ in
  \cite[Remark 9.45]{JLR}. 
While we do not allow loops in our current work,
for several reasons, we
  believe that the results here can be extended to multigraphs with loops.
 A simple and interesting case is when
  $G$ consists of a single vertex with $d/2$ loops ($d$ even).
Then $\lng$ 
  consists of the sum (union) of $d/2$ independent copies of $\mathbb
  P(n)$. Such sums have been shown to be contiguous to random
  $d$-regular (multi)graphs in \cite{SJ140}.
\end{remark}

\begin{acks}
This work was partly done at the Isaac Newton Institute, Cambridge, UK,
when CG and SJ visited in 2008, and at Institut Mittag-Leffler,
Djurs\-holm, Sweden, when SJ and AR visited in 2009.
We thank Andreas Str\"ombergsson for assistance with important 
references, and the referees for their helpful comments.
\end{acks}

\section{Notation, terminology and a summation theorem}
\label{notation}

As mentioned above, $G$ denotes a fixed connected multigraph with 
$\gv$ vertices and no loops.
For simplicity we assume that $V(G)=[\gv]\=\set{1,\dots,\gv}$.
We denote the number of edges in $G$ by $h$. 
(Often we assume $G$ to be $d$-regular, and then $h=dg/2$.) 
Let $A=A_G$ be the $\gv\times \gv$ adjacency matrix of $G$ and let
$\hA=\hA_G$ be the incidence matrix of $G$, with $\gv$ rows
and $h$ columns. Thus
\begin{equation}
  \label{haa}
\hA\hA\T = A + \DG,
\end{equation}
where $\DG$ is the
diagonal matrix with entries $d_G(i)$, $i\in V(G)$.
Denote the eigenvalues of $A$ by $\alpha_1,\ldots, \alpha_g$.

In Section~\ref{sm} we also need a directed incidence matrix
for $G$.
Give each edge in $G$ an (arbitrary) direction, and let $\aaxg$ be
the corresponding directed incidence matrix. In other words, $\aaxg$
is the $g\times h$ matrix obtained from $\hA$ by changing the sign
of one of the two 1's in each column.   Then
  \begin{equation}\label{haax}
    \aaxg\aaxg^{\, T}=D_G-A.
\end{equation}

Our version of Laplace's method (Theorem~\ref{TA2}) involves
lattices.  
A \emph{lattice} is a discrete subgroup of $\bbrn$. (Discrete means that the 
intersection with any bounded set in $\bbrn$ is finite.) It is well-known that 
every lattice $\cL$ is isomorphic (as a group) to $\bbZ^r$ for some $r$ with 
$0\le r\le n$.  
The integer $r$ is called the \emph{rank} of $\cL$ and is 
denoted by $\rank(\cL)$. 
In other words, every lattice $\cL$ has a \emph{basis}, \ie{} a sequence
$x_1,\dots,x_r$ of elements of $\cL$ such that every element of $\cL$ has a 
unique 
representation $\sum_{i=1}^r n_ix_i$ with $n_i\in \bbZ$. Furthermore, the basis 
elements $x_1,\dots,x_r$ are linearly independent (over $\bbR$); thus the rank equals 
the dimension of the linear subspace spanned by $\cL$.

The basis is not unique (except in the trivial case $r=0$); if 
$\Xi=(\xi_{ij})$ is any 
$r\times r$ integer matrix such that the determinant $\det(\Xi)=\pm1$ (which is 
equivalent to the condition that both $\Xi$ and $\Xi\qw$ 
are integer matrices) and 
$(x_i)_1^r$ is a basis of $\cL$, then $y_i=\sum_j \xi_{ij} x_j$ 
defines another basis 
$y_1,\dots,y_r$; conversely, given $(x_i)_1^r$, every basis of $\cL$ is
obtained in this way by some such matrix $\Xi$.

A \emph{unit cell} of the lattice $\cL$ is the set \set{\sum_{i=1}^r t_ix_i: 0\le t_i<1} for some
basis $(x_i)_i$ of $\cL$. If $\cL\subset\bbR^N$ has full rank $N$, and 
$U$ is any unit cell of
$\cL$, then $\set{x+U}_{x\in\cL}$ is a partition of $\bbR^N$.

The unit cells of a lattice $\cL$ all have the same $r$-dimensional volume (Hausdorff measure),
where $r=\rank(\cL)$; this volume is the \emph{determinant} (or \emph{covolume}) of $\cL$, and is
denoted by $\Vol(\cL)$.

If $(x_i)_{i=1}^r$ is a sequence of vectors in $\bbR^N$, the symmetric  matrix $\Gram xr$ of their
inner products is called their \emph{Gram matrix}. It is well-known that $x_1,\dots,x_r$ are
linearly independent if and only if the Gram matrix is non-singular, \ie, if and only if the
\emph{Gram determinant} $\detGram xr\neq0$.

The following results are well-known.

\begin{lemma}\label{L1}
If $(x_i)_{i=1}^r$ is a basis of a lattice $\cL$ in $\bbR^N$, then
  \begin{equation}
    \label{lg}
\detGram xr=\Vol(\cL)^2.
  \end{equation}\end{lemma}

\begin{lemma}
  \label{L5}
If $\cL_1\subseteq \cL_2$ are two lattices of the same rank, then $\cL_2/\cL_1$ is a finite group
of order $\Vol(\cL_1)/\Vol(\cL_2)$.
\end{lemma}

The \emph{Hessian} or second derivative $D^2\phi(x_0)$ of a function $\phi$ 
at a point $x_0\in\bbR^N$ is an $N\times N$ matrix; it is also naturally 
regarded as a bilinear form on $\bbR^N$. In general, if $B$ is a bilinear 
form on $\bbR^N$, it corresponds to the matrix 
$\left(B(e_i,e_j)\right)_{i,j=1}^N$, where $(e_i)_{i=1}^N$ is the standard 
basis. We define the determinant $\det(B)$ as 
$\det\left(B(e_i, e_j)\right)_{i,j=1}^N$, and note that if $z_1,\dots,z_N$
is any basis in $\bbR^N$, then
\begin{equation}\label{detB}
 \det(B)=\frac{\det\left(B( z_i, z_j)\right)_{i,j=1}^N}
{\det(\innprod{ z_i, z_j})_{i,j=1}^N}.
\end{equation}

We are interested in the restriction to a subspace. If $B$ is a 
bilinear form on $\bbR^N$ and
$V\subseteq\bbR^N$ is a subspace, we let $\det(B\rv)$ denote the determinant of $B$ regarded as a
bilinear form on $V$. By \eqref{detB}, this can be computed as
\begin{equation}\label{detBV}
 \det(B\rv)=\frac{\det\left(B(z_i, z_j)\right)_{i,j=1}^r}
{\det(\innprod{ z_i, z_j})_{i,j=1}^r}.
\end{equation}
for any basis $z_1,\dots,z_r$ of $V$.

We now state our general theorem for performing
summation over a lattice using Laplace's method.

\begin{theorem}\label{TA2}
Suppose the following:
  \begin{romenumerate}
\item \label{TA2l}
$\cL\subset\bbR^N$ is a lattice with rank $r\le N$.
\item \label{TA2v}
$V\subseteq\bbR^N$ is the
  $r$-dimensional subspace spanned by $\cL$.
\item \label{TA2w}
$W=V+w$ is an affine subspace parallel to $V$, for some $w\in\bbR^N$.
\item \label{TA2d}
$\K\subset\bbR^N$ is a compact convex set with non-empty interior $\K\inre$.
\item \label{TA2phi}
$\phi:\K\to\bbR$ is a continuous function and the restriction of $\phi$ to $\K\cap W$ has a unique
maximum at some point $x_0\in \K\inre\cap W$.
\item \label{TA2hess}
$\phi$ is twice continuously differentiable in a neighbourhood of $x_0$ 
and $H:=D^2\phi(x_0)$ is its Hessian at $x_0$.
\item \label{TA2psi}
$\psi:\K_1\to\bbR$ is a continuous function on some neighbourhood 
$\K_1\subseteq \K$ of $x_0$ with $\psi(x_0)>0$.
\item \label{TA2ln}
For each positive integer $n$ there is a vector
$\ell_n\in \bbR^N$ with $\ell_n/n\in W$,
\item \label{TA2a}
For each positive integer $n$ 
there is a positive real number $b_n$ and a function
\mbox{$a_n: (\cL + \ell_n)\cap nK \to\mathbb{R}$} such that,
as \ntoo, 
\begin{align}
  a_n(\ell)&=O\bigpar{b_ne^{n\phi(\ell /n)+o(n)}}, &&
              \ell\in (\cL+ \ell_n)\cap n\K,
\label{ta2a} \intertext{and}
  a_n(\ell)&=b_n\bigpar{\psi(\ell /n)+o(1)}{e^{n\phi(\ell /n)}},
    && \ell\in (\cL + \ell_n) \cap n\K_1,
\notag
\end{align}
uniformly for $\ell $ in the indicated sets.
  \end{romenumerate}
Then provided $\det(-H\rv)\neq0$, as \ntoo, 
\begin{equation}\label{ta2}
 \sum_{\ell\in (\cL+\ell_n)\cap n\K} a_n(\ell )
\sim \frac{(2\pi)^{r/2}\psi(x_0)} {\Vol(\cL) \det(-H\rv)\qq} b_n
 n^{r/2}
e^{n\phi(x_0)}.
\end{equation}
\end{theorem}

We remark that Theorem~\ref{TA2} can be generalised to allow
$n$ to tend to infinity along any infinite subset $I$ of
the positive integers, with the same proof.  
(Then (viii) and (ix) need only hold for every $n\in I$.)  

\section{Expected number of perfect matchings}\label{exppm}

A \emph{fractional perfect matching} of the multigraph $G$ 
is a function $f:E(G)\to [0,1]$ such that
$$
\sum_{e\ni v} f(e) = 1 \mbox{ for all } v\in V(G).
$$
Note that every $d$-regular multigraph has a trivial fractional 
perfect matching obtained by giving each edge weight $1/d$. 
We often treat $f$ as a vector $(f(e))_{e\in E(G)}$.

First, note that if there is a \pmx{} at all in a lift $\lng$ of $G$, 
then there exists a \fpm{} $f$ of $G$ such that $nf(e)$  
is an integer 
for each $e$. Indeed, suppose that $M$ is a perfect matching of a lift of
$G$. Let $\ell_e$ be the number of edges from the fiber $F_e$ in $M$, for
each edge $e\in E(G)$. Then the function $f:E(G)\to [0,1]$ defined by 
$f(e)=\ell_e/n$ is a fractional perfect matching of $G$.
Conversely, suppose that there exists  a \fpm{} 
$z=(z_e)_e$ in $G$ such that $nz_e$ is an integer for each $e$.
We may construct an $n$-lift of $G$ that contains a perfect matching
as follows: First take $nz_e$ edges above each edge $e\in E(G)$, with
all their endpoints disjoint. This yields $n$ endpoints above each
vertex $i\in G$, so we have constructed the sets $V_i$, and a perfect
matching. Extend this perfect matching to an $n$-lift by adding
further edges between $V_i$ and $V_j$ for all edges $e=ij$.
Consequently, $L_n(G)$ has a perfect matching with positive
probability if and only if there exists a \fpm{} $z$ with $nz$ integer-valued. 
In the sequel, for a given graph $G$ we consider only those
values of $n$ for which this holds, since otherwise trivially $X_G=0$.

\begin{remark}
  It seems to be an interesting problem to characterize the set of  
  such $n$ for a given graph, but this is outside the scope of the
  present paper, and we note only the following examples:
If $G$ itself has a perfect matching 
then every $n$ is allowed. On the other hand, if $g$ is odd, then only
even $n$ are possible. 
If $G$ is of odd order and hamiltonian, then the set of allowed $n$ is
  exactly the set of positive even integers.
If $G$ is $d$-regular, then $(1/d,\dots,1/d)$ is
  a \fpm, so every multiple of $d$ is an allowed $n$ (but there might
  be others too).
The result by Linial and Rozenman~\cite{LR} implies that for a large
  class of graphs defined there, every large even $n$ is allowed.
Note finally that if $n_1$ and $n_2$ are allowed, then so is
  $n_1+n_2$. Hence the set of allowed $n$ is always infinite, unless
  it is empty, so it makes sense to talk about asymptotic results.
\end{remark}

Suppose that there exists a \fpm{} $z=(z_e)_e$ in $G$ with 
$nz$ an integer vector.  
If a \pmx{} in $\lng$ has
$\ell_e$ edges in the fiber $F_e$ over $e$, then
$\sum_{e\ni v}\ell_e = n = n\sum_{e\ni v}z_e$ for every $e$,
so $(\ell_e)_e - nz$ belongs to the lattice $\clg$
in $\bbR^{E(G)}$ 
defined by
\begin{equation*}
  \begin{split}
  \clg&\=\Bigset{(\nu_e)_e\in\bbZ^{E(G)}:\;\sum_{e\ni v}\nu_e=0
\text{ for every $v\in V(G)$}}
\\&\phantom:
=\set{\nu\in\bbZ^{E(G)}:\hA\nu=0}.
  \end{split}
\end{equation*}
(The superscript $1$ denotes the first moment.) 
Here, and elsewhere when convenient, we think of the vectors as 
column vectors although we write them as row vectors for typographical 
reasons. 
Conversely, 
if $\ell=(\ell_e)_e$ is a vector such that $\ell-nz\in\clg$, then $\ell$ 
is an integer
vector and $\sum_{e\ni v}\ell_e=\sum_{e\ni v} nz_e=n$ for every $v$.

Given such an integer vector $(\ell_e)_e\in\clg+nz$, let us compute the
expected number of \pmx{s} in $\lng$ with
$\ell_e$ edges in the fiber $F_e$.
Clearly this number is zero unless $0\leq \ell_e\leq n$ for all $e$. 
Then the endpoints of the edges in the matching may be chosen in
\begin{equation*}
\prod_{v\in
  V(G)}\frac{n!}{\prod_{e\ni v}\ell_e!}=n!^g\prod_e(\ell_e!)\qww
\end{equation*}
ways, and for each choice, there are $\elle!(n-\elle)!$
possibilities for the fiber $F_e$, with probability $1/n!$
each.  Hence, defining $K = [0,1]^{E(G)}$ we have 
\begin{equation}\label{exg}
   \E(X_G) = \sum_{\mathbf{\ell}\in (\clg+nz)\cap nK} a_n(\ell)
\end{equation}
where 
\begin{equation*}
a_n(\ell) \= n!^{\gv-h}\, \prod_{e} \frac{(n-\ell_e)!}{\ell_e!}.
\end{equation*}
(Recall that $h$ denotes the number of edges in $G$.)

We wish to evaluate the sum \eqref{exg} asymptotically by Laplace's 
method: more precisely, by applying \refT{TA2}.  
We use Stirling's formula in the following form, valid for
all $n\ge0$, where $x\vee y\=\max(x,y)$,
\begin{equation}
  \label{stirling}
\ln(n!)=n\ln n-n+\tfrac12\ln (n\vee1)+\tfrac12\ln{2\pi}+O(1/(n+1)).
\end{equation}
Let $x_e = \ell_e/n$ for all $e\in E(G)$.
Applying (\ref{stirling})
we obtain, uniformly for $\ell\in (\clg+nz) \cap nK$,
\begin{multline*}
  \ln(a_n(\ell))
 = (g-h)\ln(n!)+\sum_{e\in E(G)} \Bigpar{\ln((n-\ell_e)!)-\ln(\ell_e!)}
\\
=(g-h)\left(n(\ln(n)-1)+\nfrac12\ln(n)+\nfrac12\ln(2\pi)+O(1/n)\right)+
\sum_{e\in E(G)} (n-2\ell_e)(\ln(n)-1)
\\
\shoveright{+n\sum_{e\in E(G)} \Bigpar{(1-x_e)\ln(1-x_e)-x_e\ln(x_e)}}
\\
+
\half\sum_{e\in E(G)} \bigpar{\ln((1-x_e)\vee n\qw)-\ln(x_e\vee n\qw)}
+\sum_{e\in E(G)}O\left( \frac1{\ell_e+1}+ \frac1{n-\ell_e+1}\right).
  \end{multline*}
Since
\begin{equation*}
  \sum_{e\in E(G)} \ell_e=\thalf \sum_v \sum_{e\ni v} \ell_e
=\thalf \sum_v n
=\thalf gn,
\end{equation*}
after cancellation,
$a_n(\ell)$ can be expressed as
\begin{equation*}
a_n(\ell)
=b_n\, \psi(\ell/n)\, \exp\left(n\phi(\ell/n)\right)
\Bigpar{1+O\Bigparfrac1{\min\elle+1}+
O\Bigparfrac1{n-\max\elle+1}}
\end{equation*}
where, for $x\in\bbR^{E(G)}$,
\begin{align}
b_n &\= (2\pi n)^{(\gv-h)/2},\\
\phi(x) &\= \sum_e \bigpar{(1-x_e)\ln (1-x_e) - x_e\ln(x_e)},\label{phi}\\
\psi(x) &\= \prod_e \left(\frac{1-x_e}{x_e}\right)^{1/2},\label{psi}
\end{align}
except that if some $x_e$  or $1-x_e$ is 0, we replace it by
$1/n$  in \eqref{psi}.
This implies that $a_n(\ell)$
satisfies  condition \eqref{ta2a} of \refT{TA2} with the above
$b_n$, $\phi$, and $\psi$.
We will now check all the remaining assumptions of \refT{TA2}.
Let
\begin{equation*}
W\=\Bigset{x=(x_e)\in\bbR^{E(G)}:\sum_{e\ni v}x_e=1
\text{ for every $v\in V(G)$}}
=\bigset{x:\hA x=(1,\dots,1)}.
\end{equation*}
As is well-known, and described in \refS{appa} in detail, the sum 
\eqref{exg} is dominated by the terms where $\phi(\ell/n)$
is close to its maximum.
In order to find the maximum, we restrict ourselves to regular multigraphs, 
where the  result is simple. (The method applies to other
graphs as well, provided one can find the maximum point(s) of
$\phi$.)

\begin{lemma}
  \label{LEXmax}
Suppose that $G$ is $d$-regular, where $d\geq 3$.
Then $\phi$
defined by \eqref{phi} has a unique maximum on
$K\cap W$=\set{x\in K:\hA x=(1,\dots,1)},
attained at the point $x^0=(1/d,\dots,1/d)$.
The maximum value is
\begin{align*}
\phi(x^0) &= \frac{g}{2}\ln\left(\frac{(d-1)^{d-1}}{d^{d-2}}\right)
,\\
\intertext{and, for $\psi$ in \eqref{psi} and the
  Hessian $D^2\phi$,}
\psi(x^0) = (d-1)^{h/2},\quad & \quad
D^2\phi(x^0) = -\frac{d(d-2)}{d-1}I.
\end{align*}
\end{lemma}
\begin{proof}

We write $\phi = \nfrac{1}{2}\sum_{v\in V(G)} \phi_v$, where
\begin{equation}
  \label{philex}
\phi_v(x_e: e\ni v)= \sum_{e\ni v} \bigpar{(1-x_e)\ln (1-x_e) - x_e\ln(x_e)}.
\end{equation}
Fix a vertex $v\in V(G)$. We rename the variables $x_e$, $e\ni v$, by
$x_1,\dots, x_d$, for convenience. Since $\phi_v$ is continuous,
it has a maximum over the compact set 
$$\Sigma_d\=\Bigset{(x_i)_i\in\oi^d:\sum_1^d x_i=1}.$$
Let $x^v\in\Sigma_d$ be a maximum point of $\phi_v$. Assume first that
$x^v$ is an interior point, i.e., that $x^v\in(0,1)^d$.
Then the function $f(y)=\phi_v(x_1^v+y, x_2^v-y, x_3^v,\dots, x_d^v)$ achieves
a maximum at $y=0$. Therefore, $f'(0)=0$ and  by the chain rule,
$$\frac{\partial \phi_v(x)}{\partial x_1}(x^v)=
       \frac{\partial \phi_v(x)}{\partial x_2}(x^v).$$
By the same argument (or by the general Lagrange multiplier method),
we have that for some constant $c_v>0$
$$\frac{\partial \phi_v(x)}{\partial x_i}(x^v)=c_v,\mbox{ for }i=1,\dots,d.$$
 But
$$\frac{\partial \phi_v(x)}{\partial x_i}(x^v)=-\ln(1-x_i)-\ln x_i-2,$$
so
$$x^v_i(1-x^v_i)=\exp\{-c_v-2\}\mbox{ for all }i=1,\dots,d.$$
This implies that the $x^v_i$'s are all at the same distance from $1/2$.
That is, for some constant $c'_v\geq 0$ we have $x^v_i=1/2\pm c'_v$ for
$i=1,\dots,d$.
Since $\sum_ix^v_i=1$ and $d\ge3$,
we have to choose the minus sign for all $i$, and thus all $x^v_i$ are
equal.
Since $x^v\in\Sigma_d$
we conclude that
$x^v_i=1/d$ for $i=1,\ldots, d$.

We also have to consider the boundary of $\Sigma_d$. If, say,
$x^v_1=0$ and $0<x^v_2<1$, then $f$ above is defined for small
positive $y$ with $f'(0+)=+\infty$, so $x^v$ cannot be a maximum point
on $\Sigma_d$. 
The only remaining points are those with all $x_i\in\set{0,1}$, but
then $\phi_v(x)=0$, while $\phi_v(1/d,\dots,1/d)>0$, so these too
cannot be (global) maximum points.
Hence $x^v$ is the unique maximum point for $\phi_v$ on $\Sigma_d$.

Setting
 $x^0=(1/d,\dots,1/d)\in \bbR^g$, we have for all $x\in K\cap W$,
$$\phi(x)\le \nfrac12\sum_v\phi_v(x^v)=\phi(x^0).$$
Moreover, the inequality is strict for all $x\neq x^0$.
This proves that $x^0$ is a unique maximum point of $\phi$ in $K\cap W$.
Clearly, $x^0$ belongs to the interior of $K$. Moreover, $\phi(x^0)$ and
$\psi(x^0)$ are given by the formulas stated in Lemma \ref{LEXmax}.

Finally, the Hessian $D^2\phi(x)$ is diagonal with entries
$(1-x_e)^{-1}-x_e^{-1}$.  Hence, at $x^0$ we have
$D^2\phi(x^0) = -\frac{d(d-2)}{d-1}I$.
\end{proof}

We have verified all assumptions of \refT{TA2},
for any neigbourhood $K_1$ of $x^0$ with
$\overline{K_1}\subset K^\circ$.
To apply formula \eqref{ta2}, we still
need to compute the rank of the lattice $\clg$ and  its determinant
$\Vol(\clg)$.

\begin{lemma}\label{LLG}
  \begin{romenumerate}
\item
If\/ $G$ is non-bipartite then the lattice
$\clg$ has rank $h-\gv$ and determinant $\Vol(\clg)=\nfrac12\det(A+\DG)\qq$.
  \item
If\/ $G$ is bipartite then the lattice
$\clg$ has rank $h-\gv+1$ and determinant $\Vol(\clg)=\det(A' + \DG')\qq$,
where the matrix $A'$ (respectively, $\DG'$) is obtained by deleting the last
row and column of $A$ (respectively, $\DG$).
  \end{romenumerate}
\end{lemma}

\begin{proof}
For $v\in V(G)$ define the vector
$x^v=\bigpar{\ett{v\in e},\, e\in E(G)}$ given by the row of the
incidence matrix $\hA$ corresponding to $v$.  For convenience,
rename these vectors $x_1,\ldots, x_g$.
Then, by \eqref{haa}, the Gram matrix of $\xxg$ is
$
\hA\hA\T=A+\DG
$.
This matrix is singular if and only if there exists a non-zero vector
$y=(y_v)\in\bbR^{V(G)}$ with $y\hA=0$. This
is equivalent to $y_i=-y_j$ for every edge $ij$, and it is
easily seen that, when $G$ is connected, such a non-zero vector $y$
exists only if
$G$ is bipartite, and that if $G$ is connected and
bipartite, there is a one-dimensional space of such solutions $y$.

Consequently, in the non-bipartite case (i),
the vectors $x_1,\dots,x_g$ are
linearly independent.  We apply~\refL{L4} with $N = h$,
$m=g$ and using the vectors $x_1,\ldots, x_g$.
Let $\cL$, $\cL\ort$ and $\cL_0$
be as in \refL{L4}. Then $\clg=\cL\ort$, and thus
$\clg$ has rank $h-g$, by \refL{L4}.
Furthermore, by \refL{L1} and \eqref{haa},
$$\Vol(\cL_0)=\bigpar{\detGram xg}\qq=\det(A+\DG)\qq.$$
Moreover, $(t_v,\, v\in V(G))$ solves \eqref{qt}
if and only if $t_v\equiv -t_w\pmod1$ for every edge $vw$.
Going around an odd cycle, we see that $t_v\equiv0$ or $t_v\equiv1/2$ for every
vertex on the cycle. Since $G$ is connected, it follows that there are
exactly two solutions to \eqref{qt}: $t_v\equiv0$ for every $v$ and
$t_v\equiv1/2$ for every $v$. Hence $q=2$ in \refL{L4}, and
the result follows.

Now suppose that $G$ is bipartite.
Then the vectors $x_1,\ldots, x_{g-1}$
are linearly independent and $x_g$ can be written as a
$\{ \pm 1\}$-combination of $x_1,\ldots, x_{g-1}$,
since the sum of vectors $x^v$ over all vertices $v$
on either side of the
vertex bipartition gives the vector $(1,1,\ldots, 1)$.
We apply~\refL{L4} with $N=h$, $m=g-1$ and using the
vectors $x_1,\ldots, x_{g-1}$.
The lemma asserts that $\clg = \cL\ort$ has rank $h-g+1$, and
$$
\Vol(\cL_0) = \bigpar{\detGram x{g-1}}\qq = \det(A' + \DG')\qq.
$$
Finally, let $w\in V(G)$ correspond to $x_g$.
If $(t_v,\, v\in V(G) \setminus \{ w\})$ solves~\eqref{qt}
then $t_u=0$ for every neighbour $u$ of $w$.
In turn this implies that $t_u=0$ for every vertex $u$ at
distance 2 from $w$, and iterating this shows that $t_u=0$
for all vertices $u$ in the connected graph $G$.
Therefore $q=1$ in~\refL{L4} and the proof is complete.
  \end{proof}

\begin{example}\label{E1K4}
When  $G=K_4$,
\begin{equation*}
\det(A+\DG)
=
\begin{vmatrix}
  3&1&1&1\\1&3&1&1\\1&1&3&1\\1&1&1&3
\end{vmatrix}
=48.
\end{equation*}
  Thus \refL{LLG}(i) says that $\clg$ has rank 2 and 
  \begin{equation*}
\Vol(\clg)=\frac{\sqrt{48}}{2}=\sqrt{12}.
  \end{equation*}
\end{example}

\begin{example}
\label{E1banana}  
Let $G=K_2^3$ be the multigraph with two vertices and three parallel
edges.  Then $A+D_G = \begin{pmatrix}3&3\\3&3\end{pmatrix}$
and deleting one row and column gives the $1\times 1$ matrix $(3)$.
Hence $\clg$ has rank 2 and $\Vol(\clg) = \sqrt{3}$, using
Lemma~\ref{LLG}(ii). 
\end{example}

We are ready to apply formula \eqref{ta2} of \refT{TA2}.

\begin{theorem}
\label{expectation}
Suppose that $G$ is  $d$-regular, where $d\ge3$.
  \begin{romenumerate}
\item
If $G$ is non-bipartite then
\begin{align*}
\E X_G &\sim \frac{2(d-1)^{d\gv/4}}{\sqrt{\mathrm{det}(A+dI)}}
\, \left(\frac{d-1}{d(d-2)}\right)^{d\gv/4-\gv/2}\,
\left(\frac{(d-1)^{d-1}}{d^{d-2}}\right)^{\gv n/2}\\
 &=
 \frac{2(d-1)^{(d-1)\gv/2}}{(d(d-2))^{d\gv/4-\gv/2}\,
   \sqrt{\mathrm{det}(A+dI)}}
\left(\frac{(d-1)^{d-1}}{d^{d-2}}\right)^{\gv n/2}.
\end{align*}
  \item
If $G$ is bipartite then
\begin{align*}
\E X_G &\sim \frac{(d-1)^{d\gv/4}}{\sqrt{\mathrm{det}(A'+dI)}}
\, \left(\frac{d-1}{d(d-2)}\right)^{d\gv/4-\gv/2 +1/2}
(2\pi n)\qq
\,
\left(\frac{(d-1)^{d-1}}{d^{d-2}}\right)^{\gv n/2}\\
 &=
 \frac{(d-1)^{(d-1)\gv/2+1/2}}{(d(d-2))^{d\gv/4-\gv/2+1/2}\,
   \sqrt{\mathrm{det}(A'+dI)}}
(2\pi n)\qq
\left(\frac{(d-1)^{d-1}}{d^{d-2}}\right)^{\gv n/2}
\end{align*}
  \end{romenumerate}
where $A'$ is obtained by deleting the last row and
column of $A$.
\end{theorem}

\begin{proof}
Let $r$ be the rank of $\clg$,
 and
recall that the Hessian $H=D^2\phi(x^0)$ is diagonal and equals
$-\frac{d(d-2)}{d-1}I$ by \refL{LEXmax}.
Thus $H\rv=-\frac{d(d-2)}{d-1}I$ too, and $\det(-H\rv)=
\bigparfrac{d(d-2)}{d-1}^r$.
Hence the result follows from \eqref{exg} and Theorem~\ref{TA2},
using  Lemmas~\ref{LEXmax} and~\ref{LLG}, and the fact
that $h=dg/2$.
\end{proof}

\begin{example}
\label{K4expectation}
For $G=K_4$, $d=3$, $\gv=4$ and thus, using  \refE{E1K4},
\begin{equation*}
 \E X_G \sim  \frac{2\cdot 2^4}{3\sqrt{48}}\, \left(\frac{4}{3}\right)^{2n}
= \frac{8}{3\sqrt{3}}\, \left(\frac{4}{3}\right)^{2n}.
\end{equation*}
\end{example}

\begin{example}
\label{bananaexpectation}
For the bipartite  
multigraph $G=K_2^3$  with two vertices and three
parallel edges we have $d=3$, $g=2$ and by Example~\ref{E1banana},
$$ 
\E X_G \sim  \frac{8}{3\sqrt{3}}\, \sqrt{\pi n}\,
      \left(\frac{4}{3}\right)^n.
$$
\end{example}

\section{The second moment of $X_G$}\label{sm}

We now work towards
an asymptotic expression for the second moment of $X_G$,
using the same approach as in the previous section.
To simplify our calculations we consider only regular
multigraphs $G$ of degree at least three. 

Given a pair $(M_1,M_2)$ of perfect matchings in $\lng$, for
a vertex $i\in V(G)$ and two (possibly equal) edges $e,f\ni i$,
let $\ell\ief$ be the number of vertices in $V_i$  whose
incident edges in $M_1$ and $M_2$ lie, respectively, in the fibers
$F_e$ and $F_f$.
Form these numbers into the $gd^2$-dimensional vector
$\ell =\ell(M_1,M_2)= \bigpar{\ell\ief : i\in[g],\,e,f\ni i}$.
Let
\begin{multline*}
V^*\= \Bigl\{\bigpar{z\ief : i\in[g],\,e,f\ni i}\in\bbR^{gd^2}:
\text{for every $e\in E(G)$ with endpoints $i$ and $j$,}
     \\
 z_{iee} = z_{jee}, \qquad
\quad \sum_{f\ni i} z_{ief} = \sum_{f\ni j} z_{jef}, \qquad
\quad \sum_{f\ni i} z_{ife} = \sum_{f\ni j} z_{jfe}
\Bigr\}.
\end{multline*}
Then the vector $\ell$ belongs to the set
\begin{equation*}
  Q\=\Bigset{(z\ief)\in V^*\cap\bbZ^{gd^2}:\sum_{e,f\ni i}z\ief=n \text{ for
  } i\in[g]}.
\end{equation*}
(The three conditions in $V^* $
follow from consideration of the edges in $M_1\cap M_2$,
$M_1$ and $M_2$, respectively.)
Fix a particular vector $z$ with $nz\in Q$. 
(By our assumption that there is a perfect matching in $\lng$,
it follows that at least one such vector exists.)
Then $Q=\clv+nz$, where $\clv$ is the lattice defined by
\begin{equation*}
\clv\=\Bigset{(\nu\ief)\in V^* \cap\bbZ^{gd^2}:\sum_{e,f\ni i}
  \nu\ief=0 \text{ for
  } i\in[g]}.
\end{equation*}
(The superscript $2$ denotes the second moment.) 

Given a pair $(M_1,M_2)$ of perfect matchings and thus a vector
$\ell\in Q$, we further define, for an edge $e\in E(G)$ and an
endpoint $i$ of $e$, 
$$
s_e=s_{ie} (\ell) = \sum_{f\ni i,\, f\neq e} \ell_{ief},\quad
t_e=t_{ie}(\ell) = \sum_{f\ni i,\, f\neq e} \ell_{ife},\quad
u_e=u_{ie}(\ell) = \sum_{f,f'\ni i;\,  f,f'\neq e} \ell_{iff'};
$$
these are the numbers of edges in the fiber $F_e$ that belong to
$M_1\setminus M_2$, $M_2\setminus M_1$ and $(M_1\cup M_2)^c$,
respectively, so they do not depend on the choice of endpoint $i$ of $e$.
We have, for every edge $e$ and endpoint $i$,
$$ s_{e} + t_{e} + u_{e} + \ell_{iee} = n.
$$

We now calculate the expected number of  pairs of perfect matchings
$(M_1,M_2)$ in $\lng$ corresponding to a given nonnegative
integer vector $\ell = (\ell_{ief})\in \clv + nz$. 
First, partition each $V_i$ into $d^2$ subsets of sizes
$(\ell_{ief})_{e,f\ni i}$; this can be done in
$$ \prod_{i=1}^g \frac{n!}{\prod_{e,f\ni i} \ell_{ief}!}
         = n!^g\prod_{i=1}^g \prod_{e,f\ni i} (\ell_{ief}!)^{-1}$$
ways.
Given these partitions there are
$$ s_{e}!\, t_{e}! \, u_{e}!\, \ell_{iee}!
$$
possibilities for the fiber $F_e$ (where $i$ is an endpoint of $e$),
with probability $1/n!$ each.
Hence the expected number  of pairs $(M_1,M_2)$ of perfect matchings
in $\lng$ which correspond to the vector $\ell$ is given by
$$
a_n(\ell) = n!^{g-dg/2}\prod_{i\in [g]}
\left(\prod_{e\ni i} \left(
   \frac{s_{e}!\, t_{e}!\, u_{e}!}
                {\ell_{iee}!}\right)^{1/2}\,
   \prod_{f\ni i,\, f\neq e}
  \frac{1}{\ell_{ief}!}\right).
$$
Thus we can write 
\begin{equation}
  \label{ex2}
 \E(X_G^2) = \sum_{\ell\in (\clv + nz)\cap nK} a_n(\ell)
\end{equation}
where $K=[0,1]^{gd^2}$.
This will allow us to apply the
same arguments as used in Section~\ref{exppm}.

We now switch to continuous variables $x\in\mathbb{R}^{gd^2}$,
where $x_{ief}$ corresponds to $\ell_{ief}/n$.
Define the functions $\sigma_{ie} = \sigma_{ie}(x)$,
$\tau_{ie} = \tau_{ie}(x)$ and $\gamma_{ie}=\gamma_{ie}(x)$
to be continuous scaled analogues of $s_{ie}$, $t_{ie}$ and $u_{ie}$
respectively.  That is, 
$$\sigma_{ie} = \sum_{f\ni i,\, f\neq e} x_{ief},\qquad
   \tau_{ie} = \sum_{f\ni i,\, f\neq e} x_{ife}, \qquad
   \gamma_{ie} = \sum_{f,f'\ni i;\, f,f'\neq e} x_{iff'},
$$
so that $\sigma_{ie}(\ell/n) = s_{ie}(\ell)/n$ and so on.
Then, applying (\ref{stirling}), it follows that $a_n(\ell)$
satisfies condition (\ref{ta2a}) of Theorem~\ref{TA2}
with
\begin{align}
b_n &= (2\pi n)^{g/2+3h/2 - d^2g/2},\notag \\
\psi(x) &= \prod_{i\in [g]} \prod_{e\ni i}
\left(\frac{\sigma_{ie}\tau_{ie} \gamma_{ie}}{x_{iee}}
  \right)^{1/4}\, \prod_{f\ni i,\, f\neq e}
  x_{ief}^{-1/2},
\notag\\
\phi(x) &= \nfrac{1}{2}\sum_{i\in [g]}\sum_{e\ni i} \bigg(
  { \sigma_{ie}\ln \sigma_{ie} + \tau_{ie}\ln\tau_{ie}
    + \gamma_{ie}\ln \gamma_{ie} -
  x_{iee}\ln x_{iee}  } \notag \\ \label{phi2}
 & \hspace*{6cm} {}
  -2\sum_{f\ni i,\, f\neq e} x_{ief}\ln x_{ief} \bigg).
\end{align}
(Again, if some $x_{ief}$, $\gs_{ie}$, $\tau_{ie}$ or $\gamma_{ie}$ is
0,  then we replace it
by $1/n$  in the definition of $\psi(x)$.)

Let $W$ be the domain defined by
\begin{equation*}
  W\=\Bigset{(x\ief)\in V^*:\sum_{e,f\ni i}x\ief=1 \text{ for
  } i\in[g]}.
\end{equation*}
We conjecture that for all connected $d$-regular multigraphs $G$
with no loops,
the function $\phi$ has a unique maximum
on $K\cap W$, attained at the point
$$ x^0 = (1/d^2,\ldots, 1/d^2).$$
Unfortunately, we have been unable to prove this, and
have only been able to verify this computationally for $d=3$.
For future reference, note that
\begin{equation}
\label{values}
\psi(x^0) = \bigpar{(d-1)d^{d-2}}^{dg},\qquad
\phi(x^0) = g\ln\left(\frac{(d-1)^{d-1}}{d^{d-2}}\right).
\end{equation}

One approach to finding the maximum of $\phi$ is to mimic the proof of
Lemma~\ref{LEXmax}.  The function $\phi$ can be written as the
sum over $i=1,\ldots, g$ of functions $\phi_i$, where the
sets of variables appearing in different $\phi_i$ are disjoint.
For convenience we drop the index $i$ and rename all variables corresponding to
vertex $i$ as $x_{ef}:=x_{ief}$, and let $\sigma_e:= \sigma_{ie}$,
$\tau_e:= \tau_{ie}$, $\gamma_e:=\gamma_{ie}$.
Then
$$
\phi_i(x)= \nfrac{1}{2}\sum_{e\ni i}
\bigg\{\sigma_e\ln \sigma_e + \tau_e\ln \tau_e
    + \gamma_e\ln \gamma_e -  x_{ee}\ln x_{ee}  -
 2\!\! \sum_{f\ni i,\,  f\neq e }x_{ef}\ln x_{ef}
\bigg\}.
$$
Since $G$ is $d$-regular and $\phi_i$ depends only on the
degree of $i$ in $G$, all the functions $\phi_i$ are equivalent
under relabelling of variables.

Now define the domain 
$$
\Sigma_{d^2} = \bigg\{ (x_{ef})_{e,f\ni i} \in\oi^{d^2} :
            \sum_{e,f\ni i} x_{ef} = 1\bigg\}.
$$
It suffices to prove that $\phi_i$ 
has a unique maximum on $\Sigma_{d^2}$ 
attained at the point
$(1/d^2,\ldots, 1/d^2)$.
Applying the Lagrange multiplier method to $\Sigma_{d^2}$, we
see that at an interior maximum point, all partial derivatives of
$\phi_i$ must
be equal.  This gives ${d^2}-1$ (non-linear) equations
(together with $\sum_{e,f}x_{ef}=1$)
to be solved for $d^2$ variables.
We tried to solve this system using \maple.
Unfortunately, \maple\ seems unable to handle the computations
for $d\geq 4$.
Hence we only have the desired result for $d=3$.

\begin{lemma}
  \label{LEX2max}
If $G$ is 3-regular then
the function $\phi$ defined by
$(\ref{phi2})$
has a unique maximum on $K\cap W$ attained at the point
$(1/9,\dots, 1/9)\in\mathbb{R}^{9g}$.
\end{lemma}

\begin{proof}
As explained above, we consider only the function $\phi_i$ for a fixed
vertex $i$.  Using \maple, we solved for points in 
$\bigset{(x_{ef})_{e,f}:\sum_{e,f}x_{ef}=1}$ 
where all the 9 partial derivatives of $\phi_i$ are equal.
Exactly four solutions were found, of which
only one lies in $[0,1]^9$, giving the point
$x^0 = (1/9,\ldots, 1/9)\in\Sigma_9$. 
(The other three solutions each contain
both positive and negative entries.)
We have $\phi(x^0)=\ln(4/3)$.

It remains to consider the boundary, where one or several $x_{ef}=0$.
If $x_{ee}=0$ and $\gamma_f>0$ for $f\neq e$,  then
$\frac\partial{\partial x_{ee}} \phi(x)=+\infty$, and thus $x$ is not
a maximum point.
Similarly, $x$ cannot be a maximum point if $x_{ef}=0$, where $e\neq
f$ and at most one of $\gs_e$, $\tau_f$ and $\gamma_{f'}$ (where $f'$
is the third index) vanishes.
It is easily seen that the only remaining cases are when the only
non-zero variables (after relabelling the indices as $1,2,3$ in some
order)
are \set{x_{12},x_{21}}, \set{x_{11},x_{22},x_{33}}
or \set{x_{11},x_{12},x_{13}}, or a subset of one of these. In the
first case we have $\phi=0$. In the two latter cases, $\phi_i$ equals,
after relabelling, 
$\frac12\phi_v$ defined in
\eqref{philex} (at the corresponding step of the first moment
calculation), and thus the maximum over one of these sets is
$\frac12\ln(4/3)<\phi(x_0)$. (We omit the details.)
Hence, there is no global maximum on the boundary.

Consequently,
$x^0$ is the unique maximum point of $\phi_i$ on 
$\Sigma_{9}$. 
Arguing as in Lemma~\ref{LEXmax} completes the proof.
\end{proof}

Let $V=W- z$ be the subspace spanned by $\clv$, i.e.,
\begin{equation*}
  V\=\Bigset{(x\ief)\in V^*:\sum_{e,f\ni i}x\ief=0
\text{ for  } i\in[g]}.
\end{equation*}

\begin{theorem}
\label{conditional}
Suppose that $G$ is $d$-regular, where $d\geq 3$.
If the function $\phi$
defined in \eqref{phi2} has a unique maximum on $K\cap W$ at
$x^0 = (1/d^2,\ldots, 1/d^2)$, then
\begin{align*}
\E(X_G^2)
&\sim\frac{ \left( (d-1)d^{d-2}\right)^{dg}}{\det(\clv)\det(-H\rv)\qq}
     (2\pi n)^{r/2+g/2+3dg/4-d^2g/2}\,
     \left(\frac{(d-1)^{d-1}}{d^{d-2}}\right)^{gn},
\end{align*}
where
$r$ is the rank of $\clv$ and
$H = D^2\phi(x^0)$ is the Hessian of $\phi$ at $x^0$, provided the
determinant in the denominator is non-zero.
In particular,
this expression holds for all 3-regular connected graphs $G$.
\end{theorem}

\begin{proof}
This is now an immediate consequence of \refT{TA2}, using
\eqref{ex2} and \eqref{values}.
The final statement follows from Lemma~\ref{LEX2max}.
\end{proof}

It remains to calculate the determinants of $\clv$ and $-H\rv$, and
the rank $r$.
In the non-bipartite case, part of this is covered by the next lemma.

\begin{lemma}\label{Llattice}
Suppose that $G$ is non-bipartite and $d$-regular, where $d\geq 3$.
Recall that $h$ denotes the number of edges in $G$, so $h=dg/2$.
Then the lattice $\clv$ has rank $d^2g-(g+3h)=d^2g-g-3dg/2$ and
  determinant
  \begin{equation*}
    \begin{split}
\det(\clv)
&=2^{3h/2-3g/2-2}\bigpar{d(d-2)}^{h/2-g/2}
\det(dI+A)\det(d(2d-3)I-A)\qq
\\&
=2^{3h/2-3g/2-2}\bigpar{d(d-2)}^{h/2-g/2}
\prod_{i=1}^g(d+\ga_i)(d(2d-3)-\ga_i)\qq,
    \end{split}
  \end{equation*}
where $\ga_1,\dots,\ga_g$ are the eigenvalues of $A$.
\end{lemma}

\begin{proof}
The linear space $V$ spanned by $\clv$ is the subspace of
$\bbR^{gd^2}$ orthogonal to the following $g+3h$ vectors:
\begin{itemize}
  \item
one vector $x^{0j}$ for every $j\in V(G)$, with $x^{0j}\ief=\ett{i=j}$.
  \item
one vector $x^{1\eps}$ for every $\eps\in E(G)$, with
$x^{1\eps}\ief=\ax_{i\eps}\ett{e=f=\eps}$.
  \item
one vector $x^{2\eps}$ for every $\eps\in E(G)$, with
$x^{2\eps}\ief=\ax_{i\eps}\ett{e=\eps\neq f}$.
  \item
one vector $x^{3\eps}$ for every $\eps\in E(G)$, with
$x^{3\eps}\ief=\ax_{i\eps}\ett{e\neq\eps=f}$.
\end{itemize}
Relabel these vectors (in this order) as $x_1,\dots,x_{g+3h}$. Then
their Gram matrix $\gG$ can be written in block form, with blocks of
dimensions $g,h,h,h$:
\begin{equation*}
  \gG=
  \begin{pmatrix}
d^2I & \aax & (d-1)\aax     & (d-1)\aax \\
\aax^{\, T} & 2I & 0 & 0 \\
(d-1)\aax^{\, T} & 0 & 2(d-1)I & \aax^{\, T}\aax-2I \\
(d-1)\aax^{\, T} & 0 & \aax^{\, T}\aax-2I  & 2(d-1)I
  \end{pmatrix}.
\end{equation*}
In order to evaluate the Gram determinant $\det(\gG)$, we may make an
orthogonal change of basis in the first component $\bbR^g$, and
another orthogonal change of basis in each of the components $\bbR^h$
(we choose the same change in all three).
It is well-known that we can make such changes of basis such that
any given $g\times h$ matrix $B$ obtains the form of a
diagonal $g\times g$ matrix $D_s$ with
$h-g$ additional columns of 0's; this is known as the singular value
decomposition of $B$, and is easily seen by choosing an orthonormal
basis $z_1,\dots,z_h$ in $\bbR^h$ such that $B^TB$ is diagonal, 
and then choosing an orthonormal basis in $\bbR^g$ containing
the vectors $B z_i/\|B z_i\|$, for all $i$ such that $B z_i\neq0$.
We choose such bases for $B=\aax$. The diagonal entries $s_1\dots,s_g$
of $D_s$ can be assumed to be non-negative, and they are identified by
the fact that the eigenvalues of $BB^T=\aax\aax^{\, T}$ are \set{s_i^2}.
By \eqref{haax}, we thus have
\begin{equation}\label{singular}
  s_i^2=d-\ga_i.
\end{equation}
 Hence, with $\dsx=(D_s,\;  0)$ a $g\times h$ matrix with
non-zero elements given by \eqref{singular},
\begin{equation}\label{sop}
\det{\gG}=
  \begin{vmatrix}
d^2I & \dsx & (d-1)\dsx     & (d-1)\dsx \\
\dsx^T & 2I & 0 & 0 \\
(d-1)\dsx^T & 0 & 2(d-1)I & \dsx^T\dsx-2I \\
(d-1)\dsx^T & 0 & \dsx^T\dsx-2I  & 2(d-1)I
  \end{vmatrix}.
\end{equation}
Since $D_s$ is a diagonal matrix, we can reorder the rows and columns
in \eqref{sop} so that we obtain a block diagonal matrix with $g$
$4\times4$ blocks
\begin{equation}
\gG_i\=
  \begin{pmatrix}
d^2 & s_i & (d-1)s_i    & (d-1)s_i \\
s_i & 2 & 0 & 0 \\
(d-1)s_i & 0 & 2(d-1) & s_i^2-2 \\
(d-1)s_i & 0 & s_i^2-2  & 2(d-1)
  \end{pmatrix}
\end{equation}
and $h-g$ identical $3\times3$ blocks
\begin{equation}
\gG_0\=
  \begin{pmatrix}
 2 & 0 & 0 \\
 0 & 2(d-1) & -2 \\
 0 & -2  & 2(d-1)
  \end{pmatrix}.
\end{equation}
Hence, by straightforward calculations,
\begin{equation}\label{detgg}
  \begin{split}
\det(\gG)
&
=\det(\gG_0)^{h-g}\prod_{i=1}^g \det(\gG_i)
\\&
=\bigpar{8d(d-2)}^{h-g}
\prod_{i=1}^g (2d-s_i^2)^2\bigpar{2d^2-4d+s_i^2}
\\&
=\bigpar{8d(d-2)}^{h-g}
\prod_{i=1}^g (d+\ga_i)^2(d(2d-3)-\ga_i)
  \end{split}
\end{equation}
Since $G$ is non-bipartite, $-d<\ga_i\le d$ for every $i$, and thus
\eqref{detgg} shows that $\det(\gG)\neq0$. Hence, the vectors
$x_1,\dots,x_{g+3h}$, or in different notation 
\begin{equation}
\label{xvectors}
 \{ x^{0j} : j\in V(G)\} \cup \{ x^{1\eps} ,x^{2\eps} ,x^{3\eps} 
           : \eps \in E(G)\},
\end{equation}
are linearly independent, so they form a basis in $V^\perp$.

We apply \refL{L4}, with $N=d^2g$, $m=g+3h=g+3dg/2$, and using
the vectors $x_1,\ldots, x_{g+3h}$ in (\ref{xvectors}).
Then $\clv=\cL^\perp$. Hence, $\rank(\clv)=N-m=d^2g-g-3h$.
We have $\det(\cL_0)=\det(\gG)\qq$ by \refL{L1}. 
Finally, we claim that
there are 4 solutions (mod 1) to \eqref{qt}: 
if we let $t_{0j}$ denote the
coefficient of $x^{0j}$, and so on, the solutions
have $t_{0j}=t_0$ for all $j$ and $t_{1\eps}=t_1$, $t_{2\eps}=t_2$,
$t_{3\eps}=t_3$ for all $\eps$, where $(t_0,t_1,t_2,t_3) = (0,0,0,0)$,
$(0,0,\frac12,\frac12)$, $(\frac12,\frac12,\frac12,0)$, or
$(\frac12,\frac12,0,\frac12)$. 
(To prove this, first consider the equations in (\ref{qt}) which
correspond to variables $x_{iee}$, and use the existence of
an odd  cycle.  This gives the possible values of $t_0$ and $t_1$.
The rest of the proof follows by considering the equations in
(\ref{qt}) which correspond to variables $x_{ief}$ for a given
vertex $i$, with $e\neq f$.)

Hence $q=4$, and \refL{L4} yields
\begin{equation*}
  \det(\clv)=\det(\cL^\perp)=\det(\gG)\qq/4.
\end{equation*}
The result follows by \eqref{detgg}.
\end{proof}

\begin{example}
  \label{K4-lattice}
For $G=K_4$, we have $d=3$,  $g=4$, $h=6$, and $A$ has the eigenvalues
$3,-1,-1,-1$. Hence \refL{Llattice} yields
$\det(\clv)=2^7\,3^{5/2}\,5^{3/2}$.
\end{example}

We believe that there is a similar result for regular bipartite
graphs, but we have not explored it. (Presumably, the rank is then
$d^2g-g-3h+2$).

Unfortunately, we have not been able to find a similar general formula
for $\det(-H\rv)$ in \refT{conditional}.  
However, this quantity can be calculated directly for a particular
graph $G$, once a basis for $\clv$ is known.

\begin{example}
\label{second-momentK4}
When $G=K_4$, using \maple\ we found a basis 
$\{ z_1,\ldots, z_{14}\}$ of $V$ and then
calculated $\det(-H\rv)=2^{-22}\,3^{28}\,5^{-1}\,11^3$
using \eqref{detBV}.  
Hence by \refT{conditional} and \refE{K4-lattice},
$$ \E(X_G^2)\sim 2^{16}\, 3^{-9/2}\, 5^{-1}\,
    11^{-3/2}\, \left(\frac{4}{3}\right)^{4n}.$$
\end{example}

\begin{example}
\label{second-moment-banana} 
When $G=K_2^3$ is the multigraph with two vertices and
three parallel edges, \maple\ computations confirmed that
$\clv$ has rank 9 and gave
$\det(\clv) = 2^4 \, 3^{3/2}$ and 
$\det(-H\rv) = 2^{-16}\, 3^{18} \,5^2$.
Hence by \refT{conditional},
$$\E(X_G^2) \sim 2^{11}\, 3^{-9/2} \, 5^{-1}\, \pi n\, \left(
    \frac{4}{3}\right)^{2n}.
$$
\end{example}

\section{Short cycles in random lifts}\label{concen}

Let $Z_k$ denote the number of cycles of length $k$ in 
$\lng$, for $k\geq 2$.
(Note that $Z_2$ is zero unless there are multiple edges in $G$.) 
To apply the small subgraph conditioning method to $X_G$,
we must understand the distribution of short cycles in
random lifts, as well as their interaction with perfect
matchings.  This will enable us to verify conditions (A1) -- (A3)
of \cite[Theorem 9.12]{JLR}, with  their $Y_n$ 
given by our $X_G$ 
(the index $n$ is suppressed), and with their $X_{kn}$ given
by our $Z_k$.

To  compute the limiting distributions in (A1) and (A2) of \cite[Theorem
9.12]{JLR}, we will use the method of moments. Moreover, for (A2) we will be 
guided by \cite[Lemma 9.17 and Remark 9.18]{JLR}, which tell us that we need 
only compute asymptotically 
$$\E(X_G\, (Z_2)_{j_2}\cdots(Z_m)_{j_m})/\E X_G,$$ 
for 
integer constants $m\ge0$ and $j_2,\ldots, j_m\geq 0$.
Here $(Z)_j$ denotes the falling factorial $Z(Z-1)\cdots (Z-j+1)$.

Let $k$ be a fixed positive integer. It is more convenient to count rooted 
oriented $k$-cycles, which introduces a factor of $2k$ into the calculations.  
A $k$-cycle in $\lng$ can be then thought of as a lift of a 
\emph{non-backtracking} closed $k$-walk in $G$, which is a 
walk $i_0e_1 i_1 e_2\ldots i_{k-1}e_k$ in $G$
such that $e_j$ is an edge of $G$ with endpoints $\{i_j, i_{j+1}\}$
and $e_j\neq e_{j-1}$, for $1\leq j \leq k$.
(Here and throughout this section, 
arithmetic on indices in $k$-walks is performed modulo $k$.) 
Note that if 
$G$ is simple then any three consecutive vertices on the walk must
all be distinct.
These walks arise in various contexts 
(see for example~\cite{ABLS, AFH, HST})
and have also been
called \emph{irreducible} \cite{Friedman} and
\emph{non-backscattering} \cite{OGS}.
Denote by $w_k$ 
the number of non-backtracking closed $k$-walks in $G$, for $k\geq 2$.

The following lemma shows that condition (A1) of
\cite[Theorem 9.12]{JLR} holds.

\begin{lemma}
\label{lambdak} Let $\lambda_k = w_k/(2k)$ for all $k\geq 2$, 
where $w_k$ is the number of non-backtracking closed $k$-walks in $G$.
Then $Z_k\sim\mathrm{Po}(\lambda_k)$, jointly for all $k\geq 2$.
\end{lemma}

\begin{proof}
Fix a non-backtracking closed $k$-walk $C = i_0e_1i_1\cdots i_{k-1}e_k$ in $G$.  The 
(oriented) 
$k$-cycle $C'=f_1f_2\cdots f_k$ in $\lng$ is a lift of $C$ if $f_j\in F_{e_j}$ for 
$j=1,\dots,k$. Hence the number of possible lifts $C'$ of $C$ is $(1+o(1))n^k$,  
and each will appear in $\lng$ with probability  $(1+o(1))n^{-k}$. It follows that
$$ \E Z_k = \sum_C\sum_{C'}\P(C'\subset L_n(G))= \frac{w_k}{2k}+o(1).$$
Similar arguments hold for higher joint factorial moments, completing the proof.
\end{proof}

For the remainder of this section we restrict our attention
to $d$-regular multigraphs with $d\geq 3$.
Next we verify condition (A2) of \cite[Theorem 9.12]{JLR}
using the approach suggested in~\cite[Remark 9.18]{JLR}.

\begin{lemma}
\label{muk} 
Suppose that $G$ is $d$-regular with $d\geq 3$, and 
for $k\geq 2$,  let
$$ \mu_k = \left(1 + \left(\frac{-1}{d-1}\right)^k\right)\, \lambda_k.$$
Then for any integer $m\geq 2$ and non-negative integers $j_2,\ldots, j_m$,
 $$\frac{\E(X_G\, (Z_2)_{j_2}\cdots(Z_m)_{j_m})}{\E X_G}
         \longrightarrow \prod_{i=2}^m\mu_i^{j_i}\mbox{ as }n\to\infty.$$
\end{lemma}

\begin{proof}
For ease of notation, throughout this proof we write 
$\P(M):=\P(M\subseteq L_n(G))$, 
$\P(M,C'):=\P(M\subseteq L_n(G),C'\subseteq L_n(G))$, and so on.
First we estimate $\E(X_G\, Z_k)$. We write
$$\E(X_G\, Z_k)=\sum_M\sum_C\sum_{C'}\P(M,C')=\sum_M\P(M)\sum_C\sum_{C'} \P(C'|M),$$
where the sums extend over all possible perfect matchings $M$ in $L_n(G)$, all 
non-backtracking closed $k$-walks $C$ in $G$, and all their possible lifts $C'$, respectively. 

To calculate the inner double sum, we fix a perfect matching $M_0$ and condition 
on its presence in  $L_n(G)$. Let $C={i_0}e_1{i_1}\ldots {i_{k-1}}e_k$ be a given 
non-backtracking closed $k$-walk in $G$. For a lift $C'$ of $C$  with edges
$f_1f_2\cdots f_k$, let
\begin{equation*}\xi_j(C') = \begin{cases} 1 \text{ if $f_j\in M_0$,}\\
                         0 \text{ otherwise,}
  \end{cases} \quad \mbox{ for } 1\leq j\leq k.
\end{equation*}

To estimate  the expected number of lifts of $C$ given $M_0$, we break the 
sum over all $C'$ according to the vector $\xi(C')$:
$$
\sum_{C'}\P(C'|M_0)=\sum_{u\in \{ 0,1\}^k}\,\,\,
                   \sum_{C':\xi(C')=u}\P(C'|M_0).
$$
Let $\ell_e$ be the number of edges of $M_0$ in the fiber $F_e$, 
and say that $M_0$ is \emph{good} if
$$|\ell_e-n/d|\le n^{2/3} \mbox{ for every } e.$$ 
We may assume that $M_0$ is good, since 
the calculations for the expectation in Section~\ref{exppm}
show that the contribution from other matchings is negligible.
(Specifically, this follows from the proof of Lemma~\ref{TA1}:
in particular the fact that $S_2=o(1)$, $S_3=o(1)$, using
notation from that proof.)

Hence, for a given $u=(u_1,u_2,\ldots, u_k)\in \{ 0,1\}^k$,
$$\P(C'|M_0)\sim\left(\frac1{n-n/d}\right)^{k-\sum_iu_i}.$$
Let $t_{00}(u)$ and $t_{01}(u)$ be the numbers of substrings $00$ and $01$ in $u$, 
respectively.
Next we prove that the number of lifts $C'=f_1\dotsm f_k$ of $C$ such that $\xi(C')=u$ 
is asymptotically equal to
$$\left(n-\frac{2n}d\right)^{t_{00}(u)}\left(\frac
nd\right)^{t_{01}(u)}.
$$
Indeed, let $V_{ie}$ be the set of endpoints 
in $V_i$ of the $\ell_e$ edges in 
$M_0\cap F_e$, for $i$ incident to $e\in E(G)$. If, say, $u_1=u_2=0$, which means that 
both, $f_1$ and $f_2$, are not in $M_0$, then we can choose the end of $f_1$ in 
$V_{i_1}$ from $V_{i_1}\setminus(V_{i_1e_1}\cup V_{i_1e_2})$, and 
$|V_{i_1}\setminus(V_{i_1e_1}\cup V_{i_1e_2})|\sim n-2n/d$ since we assume that
$M_0$ is good. Similarly, if $u_1=0$ and $u_2=1$, which means that $f_1\not\in M_0$ 
but $f_2\in M_0$, then we have to choose the end of $f_1$ from $V_{i_1e_2}$, a set of 
size $\sim n/d$.  
Note also that if $u_1=1$ then we must have $u_2=0$, and if we 
have already selected the end $w$ of $f_1$ in $V_{i_0}$, then the other end of $f_1$ is 
completely determined as the partner of $w$ in $M_0$.

Multiplying these two expressions together yields that
$$\sum_{C':\xi(C')=u}\P(C'|M_0)=b_{u_1u_2}\cdots
               b_{u_{k-1}u_k}b_{u_ku_1}+o(1),$$
 where $b_{00}, b_{01}, b_{10},b_{11}$ form the matrix
\renewcommand{\arraystretch}{1.5}
\begin{equation*}  B = \begin{pmatrix} \frac{d-2}{d-1} & \frac{1}{d-1}\\
  1 & 0\end{pmatrix}.\end{equation*}
Note that $B$ has eigenvalues 1 and $-1/(d-1)$. Summing over all 
$u=(u_1,\ldots, u_k)$, we find that the conditional expected number of lifts of $C$ is
$$\sum_{C'}\P(C'|M_0)=\mathrm{Tr}(B^k)+o(1) =
  1 + \left(\frac{-1}{d-1}\right)^k+o(1).$$
Hence the expected number of $k$-cycles in $L_n(G)$, conditioned on the existence 
of a given good
perfect matching $M_0$, is asymptotically equal to 
$$\sum_C\sum_{C'} \P(C'|M_0)\sim \mu_k:=
   \left(1 + \left(\frac{-1}{d-1}\right)^k\right)\frac{w_k}{2k}=
   \left(1 + \left(\frac{-1}{d-1}\right)^k\right)\lambda_k.$$
Finally,
$$\E (X_G \, Z_k)\sim\sum_M\P(M)\mu_k=\mu_k\E X_G.$$

All the above calculations work similarly for higher factorial moments and 
yield the desired result.
\end{proof}

Denote a directed edge of $G$ by $(e,i,j)$, where $e\in E(G)$ 
is incident to $i,j\in V(G)$ and $i\neq j$; this denotes $e$ 
directed from $i$ to $j$. Now let $R$ be the $dg\times dg$ matrix
with rows and columns indexed by directed edges of $G$, and
$$ R_{(e,i,j),(f,p,q)} = \begin{cases} 1 & \text{ if $p=j$ and $f\neq e$,}\\
         0 & \text{ otherwise. }
  \end{cases}
$$
(Here $R$ is the adjacency matrix of a version of the directed line 
graph of $G$, where $U$-turns are forbidden.) Then 
\begin{equation}\label{atheta}
w_k = \mathrm{Tr}(R^k) =
             \theta_1^k + \cdots + \theta_{dg}^k,
\end{equation}
where $\theta_1,\ldots, \theta_{dg}$ are the eigenvalues of $R$. Note that 
$d-1$ is an eigenvalue of $R$ with eigenvector $(1,1,\ldots, 1)^T$; since 
$R$ has non-negative entries, this is the eigenvalue with largest modulus.
Now for $k\geq 2$, the quantity $\mu_k$ defined in Lemma~\ref{muk} equals
$$ \mu_k = (1+ \delta_k)\lambda_k \mbox{, where }
   \delta_k = \left(\frac{-1}{d-1}\right)^k > -1.$$
Therefore the quantity $\sum_k \lambda_k \delta_k^2 $ in condition 
(A3) of~\cite[Theorem 9.12]{JLR} is
\begin{align*}
\sum_k \lambda_k \delta_k^2 &= \sum_{k\geq 1} \frac{w_k}{2k\, (d-1)^{2k}}
= \sum_{k\geq 1} \frac{1}{2k} \sum_{t=1}^{dg} \left(
           \frac{\theta_t}{(d-1)^2}\right)^k\\
   &= -\frac{1}{2}\sum_{t=1}^{dg}\ln\left(1-\frac{\theta_t}{(d-1)^2}\right),
\end{align*}
which is finite as required. Furthermore,
\begin{align}
\exp\left(\sum_k \lambda_k \delta_k^2 \right)
&=(d-1)^{dg}\left(\prod_{t=1}^{dg}((d-1)^2-\theta_t)\right)^{-1/2}\notag\\
&=(d-1)^{dg}\det\bigpar{(d-1)^2I-R}\qqw.
\label{eld}
\end{align}
In order to assist with the verification of condition (A4) from
from~\cite[Theorem 9.12]{JLR}, we will rewrite this expression
in terms of the adjacency matrix $A$ of $G$. 
The following result was proved by Friedman~\cite{Friedman}.

\begin{lemma}
\label{Meigvals}
{\bf \cite[Theorem 10.3]{Friedman}}\
Suppose that $G$ is $d$-regular with $d\geq 3$ 
and let $\alpha_1,\ldots, \alpha_g$ be the eigenvalues of 
the adjacency matrix of $G$.
For $i=1,\ldots, g$  denote the roots of the quadratic 
$x^2-\ga_ix+d-1=0$ by $\beta_i^+$ and $\beta_i^-$. 
That is, 
$$
  \beta_i^+ = \nfrac12\alpha_i + \sqrt{\nfrac14\alpha_i^2-(d-1)},
\quad
  \beta_i^- = \nfrac12\alpha_i - \sqrt{\nfrac14\alpha_i^2-(d-1)}.
$$
Then
the eigenvalues of $R$ are $\beta_i^+$, $\beta_i^-$ for $i=1,\ldots, g$, 
together with $1$ and $-1$, the latter two repeated $g(d-2)/2$ times each.
Hence, for $k\geq 2$, the number of non-backtracking closed $k$-walks
in $G$ is given by
\begin{equation*}
w_k  = \tfrac{1}{2} g(d-2)\left( 1 + (-1)^k \right) +
        \sum_{i=1}^g \left( (\beta_i^+)^k + (\beta_i^-)^k\right).
\end{equation*}
\end{lemma}

Note that there may be repetitions among $\gb^+_i,\gb_i^-$, and some 
of these may coincide with $\pm1$. Hence 
the multiplicities of these eigenvalues may not be exactly 1 or $g(d-2)/2$: 
see~\refE{K4-magic} below.

We now use Lemma~\ref{Meigvals} to
rewrite (\ref{eld}) in terms of the
eigenvalues of the adjacency matrix of $G$.

\begin{corollary}  \label{Cw}
Suppose that $G$ is $d$-regular, with $d\geq 3$. 
The expression in $(\ref{eld})$ can be written as
\begin{align*}
 & \exp\left( \sum_{k} \lambda_k \delta_k^2\right) \\
&=
   (d-1)^{dg-g/2}((d-1)^4-1)^{-(d-2)g/4}
       \,\, \det((d-1)^3 + 1)I - (d-1)A)^{-1/2}\\
&=
   (d-1)^{dg-g/2}((d-1)^4-1)^{-(d-2)g/4}
       \,\, \prod_{i=1}^g\left((d-1)^3 + 1 - (d-1)\ga_i\right)^{-1/2}.
\end{align*}
\end{corollary}

\begin{proof}
It follows from Lemma~\ref{Meigvals} that the 
characteristic polynomial of $R$ is given by
\begin{equation*}
  \begin{split}
  \det(\gl I-R)
&=\prod_{i=1}^{dg}(\gl-\theta_i)
=(\gl-1)^{(d-2)g/2}(\gl+1)^{(d-2)g/2}\prod_{i=1}^g(\gl-\gb_i^+)(\gl-\gb_i^-)
\\&
=(\gl^2-1)^{(d-2)g/2}\prod_{i=1}^g(\gl^2-\ga_i\gl+d-1)\\
&= (\gl^2-1)^{(d-2)g/2}\, \det((\lambda^2+d-1)I - \lambda A). 
  \end{split}
\end{equation*}
The proof is completed by substituting this into (\ref{eld})
with $\lambda = (d-1)^2$. 
\end{proof}

\begin{example}
\label{K4-magic} 
When $G=K_4$ the eigenvalues of $A$ are
$\alpha_1=3,\alpha_2=\alpha_3=\alpha_4=-1$. 
By \refL{Meigvals}, the eigenvalues of $R$ are 2, 1 (three times), 
$-1$ (twice), 
and $\frac12(-1\pm \sqrt7i)$ (three times each),
so the number of non-backtracking closed $k$-walks in $K_4$ is
\begin{equation*}
w_k  = 2^k + 3 + 2(-1)^k +3\Bigparfrac{-1+\sqrt7i}{2}^k +
    3\Bigparfrac{-1-\sqrt7i}{2}^k.
\end{equation*}
Furthermore, by \refC{Cw},
$$\exp\lrpar{\sum_{k}\lambda_k\delta_k^2}
= 2^{10}\, 15^{-1}\det(9I-2A)\qqw
     = 2^{10} \, 3^{-3/2} \, 5^{-1} \, 11^{-3/2}.
$$
\end{example}

\begin{example}
\label{banana-magic}
The multigraph with two vertices connected by $d$ 
parallel edges has adjacency matrix
\begin{equation*}  A = \begin{pmatrix} 0 & d\\
  d & 0\end{pmatrix}.\end{equation*}
We have $\gb_1^\pm,\gb_2^\pm=\pm(d-1),\,\pm1$ and by \refL{Meigvals}, 
the matrix
$R$ has eigenvalues $\pm(d-1)$ and $\pm1$, the latter with
mulitiplicities $d-1$.
Hence $w_k=2(d-1)^k+2(d-1)$ if $k\ge2$ is even, and
$w_k=0$ if $k$ is odd.
\refC{Cw} yields, after some algebra, 
$$\exp\left(\sum_k \lambda_k \delta_k^2 \right)=
     (d-1)^{2d-1}d^{-d/2}(d-2)^{-d/2}(d^2-2d+2)^{-d/2+1/2}.$$
For example, when $d=3$ this is $2^5 3^{-3/2} 5^{-1}$, while for $d=4$ it is 
$2^{-{15/2}} 3^7 5^{-3/2}$. 
\end{example}

To complete this section, we prove a concentration result for 
the number of perfect matchings in $\lng$ when $G=K_4$ and 
when $G$ is the multigraph $K_2^3$ with 2 vertices and 3 parallel edges.
We conjecture that the analogous result is true for any connected
$d$-regular multigraph $G$ with no loops, where $d\geq 3$, 
with $\delta_k = -(1/(d-1))^k$.

\begin{corollary}\label{mainK4}
For $k\geq 3$ let $w_k$ be the number of non-backtracking
closed walks of length 
$k$ in $K_4$, and define $\lambda_k = w_k/2k$. Further, let 
$Y_k$ be a Poisson random variable with
expectation $\lambda_k$, with $\set{Y_k}_k$ independent, 
and define $\delta_k = (-1/2)^k$.  Then 
with $G = K_4$,
$$\frac{X_G}{\E X_G}\dto
    W:=\prod_{i=3}^{\infty}
          \left(1+\delta_i\right)^{Y_i}e^{-\lambda_i\delta_i}. 
$$ 
\end{corollary}

\begin{proof}
Let $X = X_{K_4}$.
It follows from 
Examples \ref{K4expectation} and~\ref{second-momentK4} that 
$$ \frac{\E(X^2)}{(\E X)^2} 
   \sim 2^{10}\, 3^{-3/2}\, 5^{-1}\, 11^{-3/2}.
$$
By comparing with Example~\ref{K4-magic}, we find that (A4) 
of~\cite[Theorem 9.12]{JLR} is satisfied:  that is, 
$$ \frac{\E X^2}{(\E X)^2} \to
             \exp\left( \sum_k \lambda_k \delta_k^2\right) \quad
     \mbox{ as $n\to\infty$}.
$$
The other conditions of \cite[Theorem 9.12]{JLR} hold, as follows
from Lemmas~\ref{lambdak} and \ref{muk}. 
Applying \cite[Theorem 9.12]{JLR}  completes the proof.
\end{proof}

The same argument applies for the multigraph with two
vertices and three parallel edges, this time using
Examples~\ref{bananaexpectation}, 
\ref{second-moment-banana} and \ref{banana-magic}, leading to
the following.  

\begin{corollary}
\label{mainbanana}
Recall that $K_2^3$ denotes the multigraph with two vertices
and three parallel edges.
For $k\geq 2$  
let $w_k$ be the number of non-backtracking closed walks of length $k$, 
and define $\lambda_k = w_k/2k$. Further, let 
$Y_k$ be a Poisson random variable with
expectation $\lambda_k$, with $\set{Y_k}_k$ independent,  
and define $\delta_k = (-1/2)^k$.  Then 
with $G = K_2^3$,  
$$\frac{X_G}{\E X_G}\dto
    W:=\prod_{i=1}^{\infty}
          \left(1+\delta_{2i}\right)^{Y_{2i}}e^{-\lambda_{2i}\delta_{2i}}.
$$ 
\end{corollary}

It is immediate that the limiting distribution $W$ satisfies
$W > 0$ (with probability 1) in both
Corollary~\ref{mainK4} and \ref{mainbanana}. 
Hence $\lng$ a.a.s.\ has a perfect matching, for both $G=K_4$
and $G=K_2^3$.
This also follows from~\cite{LR}.

\section{Summation by Laplace's method}\label{appa}

In this section we prove our main approximation tool, Theorem \ref{TA2},
which performs a summation over lattice points. 
We will require a little more theory about lattices.
The following surprising duality was 
proved by McMullen \cite{McM}.  (See also \cite{Schnell}.)

\begin{lemma}  \label{L4a}
Let $V$ be a subspace of $\bbrn$ and let $V\ort$ be its orthogonal complement. 
Let $\cL$ and $\cL\ort$ be the lattices $V\cap\bbzn$ and $V\ort\cap\bbzn$, 
and assume that the rank of $\cL$ equals the dimension of $V$ 
(\ie, that $\cL$ spans $V$). Then $\cL\ort$ has rank
$\dim(V\ort)=N-\dim(V)$ and
\begin{equation*}
  \Vol(\cL\ort)=\Vol(\cL).
\end{equation*}
\end{lemma}
For our purposes we need a simple extension.

\begin{lemma}  \label{L4}
  Let $0\le m\le N$.
Let $\xxm$ be linearly independent vectors in $\bbZ^N$. Let $V$ be the subspace of $\bbrn$ spanned
by $\xxm$ and let $V\ort$ be its orthogonal complement; thus
\begin{equation*}
  V\ort=\set{y\in\bbrn:\innprod{y,x_i}=0\text{ for }i=1,\dots,m}.
\end{equation*}
Let $\cL$ and $\cL\ort$ be the lattices $V\cap\bbzn$ and $V\ort\cap\bbzn$, and let $\cL_0$ be the
lattice spanned by $\xxm$ (\ie, the set 
\set{\sum_{i=1}^m n_ix_i:n_i\in\bbZ} of integer combinations).
Then $\cL\ort$ has rank $N-m$ and
\begin{equation*}
  \Vol(\cL\ort)=\Vol(\cL)=\Vol(\cL_0)/q,
\end{equation*}
where $q$ is the order of the finite group $\cL/\cL_0$. 
Explicitly, $q$ is the number of solutions
$(\ttm)$ in $(\bbR/\bbZ)^m$ (or $(\bbQ/\bbZ)^m$) of the system 
\begin{equation}\label{qt}
  \sum_i x_{ij} t_i \equiv 0 \pmod1,
\qquad j=1,\dots,N,
\end{equation}
where $x_i = (x_{ij})_{j=1}^N$ for $i=1,\ldots, m$.
\end{lemma}

\begin{proof}
Since $\rank(\cL)=m=\dim(V)$, we can apply \refL{L4a} and conclude that 
$\rank(\cL\ort)=N-m$ and $\Vol(\cL\ort)=\Vol(\cL)$.

Next, $\cL_0\subseteq V\cap\bbzn=\cL$; moreover, $\cL_0$ and $\cL$ both span $V$ and have thus the
same rank. Hence \refL{L5} shows that $\cL/\cL_0$ is finite and $\Vol(\cL)=\Vol(\cL_0)/q$. Note
further that $\cL\subseteq V = \set{\sum_i t_i x_i :t_i\in\bbR}$ and thus
$$
q=|\cL/\cL_0| = \Bigabs{\Bigset{(t_i)\in[0,1)^m:\sum_i t_ix_i\in \cL}}.
$$
Furthermore,
$$
\sum_i t_ix_i\in \cL \iff \sum_i t_ix_i\in \bbzn \iff \sum_i x_{ij} t_i \equiv 0 \pmod1 \text{ for
}j=1,\dots,J,
$$
and the characterization of $q$ follows.
\end{proof}

The proof of Theorem \ref{TA2} involves reduction to a special case,
 which we prove first.
\begin{lemma}\label{TA1}
Suppose the following: 
  \begin{romenumerate}
\item \label{TA1l}
$\cL\subset\bbR^r$ is a lattice with full rank $r$.
\item \label{TA1d}
$\K\subset\bbR^r$ is a compact convex set with non-empty interior $\K\inre$.
\item \label{TA1phi}
$\phi:\K\to\bbR$ is a continuous function with a unique maximum at some interior point $x_0\in
\K\inre$.
\item \label{TA1hess}
$\phi$ is twice continuously differentiable in a neighbourhood of $x_0$ and the Hessian
$H:=D^2\phi(x_0)$ is strictly negative definite.
\item \label{TA1psi}
$\psi:\K_1\to\bbR$ is a continuous function on some neighbourhood $\K_1\subseteq \K$ of $x_0$ with
$\psi(x_0)>0$.
\item \label{TA1ln}
For each positive integer $n$ there is a vector $\ell_n \in \bbR^r$.
\item \label{TA1a}
For each positive integer $n$ there is a positive real number
$b_n$ and a function \mbox{$a_n: \llnx \cap nK \to\mathbb{R}$}
such that, as \ntoo, 
\begin{align}
  a_n(\ell)&=O\bigpar{b_ne^{n\phi(\ell /n)+o(n)}},
 && \ell\in \llnx\cap n\K,
\label{t1a} \intertext{and}
  a_n(\ell)&=b_n\bigpar{\psi(\ell /n)+o(1)}{e^{n\phi(\ell /n)}}, &&
\ell\in \llnx\cap n\K_1, \label{t1b}
\end{align}
uniformly for $\ell $ in the indicated sets.
  \end{romenumerate}
Then, as \ntoo,
\begin{equation}\label{t1}
 \sum_{\ell\in (\cL+\ell_n)\cap n\K} a_n(\ell )
\sim \frac{ (2\pi)^{r/2}\psi(x_0)} {\Vol(\cL) \det\bigpar{-H}\qq} \, b_n n^{r/2} e^{n\phi(x_0)}.
\end{equation}
\end{lemma}

\begin{proof}
We begin with a few simplifications. We may obviously assume 
that $b_n=1$.  Furthermore,
  by subtracting $\phi(x_0)$ from $\phi$, and dividing $a_n(\ell)$ by
  $e^{n\phi(x_0)}$, we may suppose that $\phi(x_0)=0$.

Since $x_0$ is an interior maximum point, the gradient $D\phi(x_0)$ vanishes, 
and a Taylor
expansion at $x_0$ shows that, using \ref{TA1hess}, as $|x-x_0|\to0$,
\begin{align}\label{taylor}
  \phi(x)
&=\thalf \innprod{x-x_0,D^2\phi(x_0)(x-x_0)} + o(|x-x_0|^2)
\\&
\le -\cc|x-x_0|^2+ o(|x-x_0|^2) \notag
\end{align}
for some positive constant $\ccx$.
Consequently, there exists $\gd>0$ such that the neighbourhood 
$\set{x:|x-x_0|\le\gd}$ is contained in $\K_1$ and
\begin{equation}\label{a1}
  \phi(x)\le -\cc|x-x_0|^2, \ccdef\ccai
\qquad |x-x_0|<\gd
\end{equation}
for some positive constant $\ccx$.
We divide the sum in \eqref{t1} into three parts:
\begin{align*}
  S_1:=\hskip-1em\sum_{|\ell /n-x_0|< n^{-1/3}},
&&
  S_2:=\hskip-1em\sum_{n^{-1/3}\le|\ell /n-x_0|< \gd},
&&
  S_3:=\hskip-4pt\sum_{|\ell /n-x_0|\ge\gd}.
\end{align*}
In the sum $S_2$ we use \eqref{t1b} and \eqref{a1}; thus each term is
\begin{equation*}
 a_n(\ell)=O(e^{n\phi(\ell/n)})
=O(e^{-\ccx n\qqq}).
\end{equation*}
Since the number of terms is $O(n^r)$, we obtain $S_2=o(1)$.

Similarly, by compactness, if $|x-x_0|\ge\gd$, then $\phi(x)\le -\cc$
for some positive constant $\ccx$. Consequently, for large $n$,
\eqref{t1a} shows that each term in $S_3$ is
\begin{equation*}
 a_n(\ell)=O(e^{n\phi(\ell/n)+\ccx n/2})
=O(e^{-\ccx n/2}).
\end{equation*}
Again, the number of terms is $O(n^r)$ and we obtain $S_3=o(1)$.

We convert the sum $S_1$ into an integral by picking a unit cell 
$U$ of the lattice $\cL$ and
defining $a_n(y)\=a_n(\ell )$ for $y\in U+\ell $, $\ell\in \cL+\ell_n$. Let $Q_n:=\bigcup_{|\ell
/n-x_0|< n^{-1/3}} (U+\ell )$, and let $\tqn:=\set{z:nx_0+\sqrtn z\in Q_n}$. Then
\begin{equation}\label{s1}
  S_1
=\Vol(\cL)\qw \int_{Q_n} a_n(y) \dd y =\Vol(\cL)\qw n^{r/2} \int_{\tqn} a_n\bigpar{nx_0+\sqrtn z}
\dd z.
\end{equation}
Note that $Q_n$ is roughly a ball of radius $n^{2/3}$ centered at $nx_0$, and $\tqn$ is roughly a
ball of radius $n^{1/6}$ centered at $0$.

If $y\in Q_n$, then $|y/n-x_0|\le n^{-1/3}+O(n\qw)$. Since the gradient $D\phi(x_0)=0$,
\ref{TA1hess} implies that for $x\in Q_n/n$,
\begin{equation}\label{d1}
  |D\phi(x)| = O(|x-x_0|)= O(n\qqqw).
\end{equation}
If $y\in U+\ell \subset Q_n$, then $|y/n-\ell /n|=O(1/n)$ and \eqref{d1} implies
\begin{equation*}
  n\phi(y/n)-n\phi(\ell /n)= O\bigpar{n n\qqqw n\qw}= O\bigpar{n\qqqw},
\end{equation*}
and thus \eqref{t1b} implies, uniformly for $y\in Q_n$,
\begin{equation}\label{a2}
  a_n(y)=a_n(\ell)=\bigpar{\psi(y/n)+o(1)}{e^{n\phi(y/n)}}.
\end{equation}

For every fixed $z\in\bbR^r$, this and the Taylor expansion \eqref{taylor} show that, as \ntoo,
using the continuity of $\psi$,
\begin{equation*}
  a_n(nx_0+\sqrtn z)\to \psi(x_0) e^{\half \innprod{z,D^2\phi(x_0)z} }.
\end{equation*}
Moreover, \eqref{a1} and \eqref{a2} provide a uniform bound, for all $z\in\bbR^r$,
\begin{equation*}
  |a_n(nx_0+\sqrtn z)\etta_{\tqn}(z)|\le \CC e^{-\ccai |z|^2}.
\end{equation*}
Further, $\etta_{\tqn}(z)\to1$ for every $z$. Hence, 
dominated convergence shows that
\begin{equation*}
  \begin{split}
 \int_{\tqn} a_n(nx_0+\sqrtn z)\dd z
&\to \int_{\bbR^r} \psi(x_0) e^{\half \innprod{z,D^2\phi(x_0)z} } \dd
 z
\\&
=\psi(x_0) (2\pi)^{r/2} \det\bigpar{-D^2\phi(x_0)}\qqw.
  \end{split}
\end{equation*}
The result follows from this and \eqref{s1}, together with the estimates $S_2=o(1)$ and $S_3=o(1)$
above.
\end{proof}

\begin{proof}[Proof of Theorem~\ref{TA2}]
First, replacing $\K$ by $\K-w$, $a_n(\ell)$ by
 $a'_n(\ell)\=a_n(\ell+nw)$, $\ell_n$ by $\ell_n-nw$,
 and translating $\phi$ and $\psi$, we reduce to the case $w=0$ and
 thus $W=V$ and $\ell_n\in V$.

Choose a lattice basis  \set{z_1,\dots,z_r}  of $\cL$. Consider the 
mapping $T:\bbR^r\to V\subseteq\bbR^N$ given by 
$(y_1,\dots,y_r)\mapsto \sum_{i=1}^r y_iz_i$, which thus maps $\bbZ^r$
onto $\cL$. We apply \refL{TA1} to $\cL'\=\bbZ^r$, $\K'\=T\qw(\K)$, 
$\phi\circ T$, $\psi\circ T$, $\ell'_n\=T\qw(\ell_n)$, and $a_n(T(k))$, 
$k\in (\cL'+\ell'_n)\cap n \K'$. The Hessian 
$D^2 (\phi\circ T)(T\qw x_0)$ equals 
$\left(H(z_i,z_j)\right)_{i,j=1}^r$, and its negative has
determinant, by \eqref{detBV} and \eqref{lg},
\begin{equation}
  \label{hess}
\det\left(- H( z_i,z_j)\right)_{i,j=1}^r 
=\det(-H\rv)\det(\innprod{z_i,z_j})_{i,j=1}^r
=\det(-H\rv)\det(\cL)^2.
\end{equation}
Hence, \eqref{ta2} follows from \refL{TA1}. 
Note that the Hessian $D^2 (\phi\circ T)(T\qw x_0)$ is
always negative semi-definite, because $x_0$ is a maximum point. 
Hence, it is negative definite unless its determinant is zero, 
which is ruled out by \eqref{hess} and the assumption
that $\det(-H\rv)\neq0$.
\end{proof}

\end{document}